\documentclass[final]{siamltex}
\usepackage{amsmath}
\usepackage{amsfonts}
\usepackage{graphicx}
\usepackage{subfigure}

\def\sec#1{Section~\ref{#1}}

\def\bF{\mathbf F}
\def\bE{\mathbf E}
\def\bH{\mathbf H}

\def\bS{\mathbf S}

\def\bx{\mathbf x}
\def\by{\mathbf y}

\def\bk{\mathbf k}

\def\b0{\mathbf 0}
\def\bTheta{\mathbf \Theta}

\title{Transverse electric scattering on inhomogeneous objects: 
singular integral equation, symbol of the operator, and matrix elements}

\author{Grigorios~P.~Zouros\thanks{School of Electrical and Computer 
        Engineering, National Technical University of Athens, 
        Athens 15773, Greece ({\tt zouros@mail.ntua.gr}).}
        \and Neil~V.~Budko\thanks{Numerical Analysis, DIAM, 
        Faculty of Electrical Engineering, Mathematics, 
        and Computer Science, Delft University of Technology, 2628 CD Delft, 
        The Netherlands ({\tt n.v.budko@tudelft.nl}).}}

\begin{document}

\maketitle

\begin{abstract}
This is a companion report for the paper ``Transverse electric scattering on 
inhomogeneous objects: Spectrum of integral operator and preconditioning'' by
the present 
authors \cite{ZourosBudko_2011}.
In this report we formulate the two-dimensional transverse electric scattering problem
as a standard singular integral equation, derive the symbol of the integral
operator for H{\"o}lder-continuous contrasts, and calculate the elements 
of the system matrix obtained after discretization via the mid-point rule.
\end{abstract}



\pagestyle{myheadings}
\thispagestyle{plain}
\markboth{G.~P.~ZOUROS AND N.~V.~BUDKO}{TE INTEGRAL EQUATION, SYMBOL, MATRIX}

\section{Derivation of the TE singular integral equation}
Although some of the material presented here may be new, in particular, the explicit expression for the operator symbol,
the methods and techniques are standard and can be found in \cite{Mikhlin, Martin_1998, Ilinski_2000}.

Starting from the frequency-domain Maxwell's equations \cite{ZourosBudko_2011} with the time convention
$e^{-i\omega t}$ we arrive at the following set of equations for the scattered fields with the
induced currents as sources:
\begin{align}
\label{eq:MatrixTE}
\begin{bmatrix}
-i\omega\varepsilon_{\rm b} & 0 & -\partial_{2}\\
0 & -i\omega\varepsilon_{\rm b} & \partial_{1} \\
-\partial_{2} & \partial_{1} & -i\omega\mu_{\rm b}
\end{bmatrix}
\begin{bmatrix}
E_{1}^{\rm sc}\\
E_{2}^{\rm sc}\\
H_{3}^{\rm sc}
\end{bmatrix}
&=
\begin{bmatrix}
-J_{1}^{\rm ind}\\
-J_{2}^{\rm ind}\\
-K_{3}^{\rm ind}
\end{bmatrix},
\\
\label{eq:MatrixTM}
\begin{bmatrix}
-i\omega\mu_{\rm b} & 0 & \partial_{2}\\
0 & -i\omega\mu_{\rm b} & -\partial_{1} \\
\partial_{2} & -\partial_{1} & -i\omega\varepsilon_{\rm b}
\end{bmatrix}
\begin{bmatrix}
H_{1}^{\rm sc}\\
H_{2}^{\rm sc}\\
E_{3}^{\rm sc}
\end{bmatrix}
&=
\begin{bmatrix}
-K_{1}^{\rm ind}\\
-K_{2}^{\rm ind}\\
-J_{3}^{\rm ind}
\end{bmatrix},
\end{align}
where
\begin{align}
\label{eq:ContrastCurrents}
\begin{split}
J^{\rm
ind}_{k}(\bx,\omega)&=-i\omega\left[\varepsilon(\bx,\omega)-\varepsilon_{\rm
b}\right]E_{k}(\bx,\omega), \;\;\;k=1,2,3;
\\
K^{\rm ind}_{m}(\bx,\omega)&=-i\omega\left[\mu(\bx,\omega)-\mu_{\rm
b}\right]H_{m}(\bx,\omega), \;\;\;m=1,2,3.
\end{split}
\end{align}
Equation \eqref{eq:MatrixTE} describes the TE case while \eqref{eq:MatrixTM} the
TM case and is the dual of \eqref{eq:MatrixTE}. The two-dimensional Fourier transform
of \eqref{eq:MatrixTE} and \eqref{eq:MatrixTM} with respect to coordinates 
$x_{1}$ and $x_{2}$ takes us from the $(\bx,\omega)$ domain to the $(\bk,\omega)$ domain,
whereby $\partial_n\rightarrow-ik_n,\,n=1,2$. In this way we arrive at the following linear
algebraic problem of the form $\mathbb A\tilde\bF=-\tilde\bS$ for the TE case:
\begin{align}
\label{ks_TE}
\begin{bmatrix}
-i\omega\varepsilon_{\rm b} & 0 & ik_2\\
0 & -i\omega\varepsilon_{\rm b} & -ik_1 \\
ik_2 & -ik_1 & -i\omega\mu_{\rm b}
\end{bmatrix}
\begin{bmatrix}
\tilde E_{1}^{\rm sc}\\
\tilde E_{2}^{\rm sc}\\
\tilde H_{3}^{\rm sc}
\end{bmatrix}
&=-
\begin{bmatrix}
\tilde J_{1}^{\rm ind}\\
\tilde J_{2}^{\rm ind}\\
\tilde K_{3}^{\rm ind}
\end{bmatrix}.
\end{align}
Solving \eqref{ks_TE} via matrix inversion, we obtain the scattered fields in
terms of the induced currents:
\begin{align}
\label{ks_TE_2}
\begin{bmatrix}
\tilde E_{1}^{\rm sc}\\
\tilde E_{2}^{\rm sc}\\
\tilde H_{3}^{\rm sc}
\end{bmatrix}=
\begin{bmatrix}
\frac{k_1^2-\omega^2\varepsilon_{\rm b}\mu_{\rm b}}{i\omega\varepsilon_{\rm b}}
& \frac{k_1k_2}{i\omega\varepsilon_{\rm b}} & ik_2\\
\frac{k_1k_2}{i\omega\varepsilon_{\rm b}} & \frac{k_2^2-\omega^2\varepsilon_{\rm
b}\mu_{\rm b}}{i\omega\varepsilon_{\rm b}} & -ik_1 \\
ik_2 & -ik_1 & i\omega\varepsilon_{\rm b}
\end{bmatrix}
\begin{bmatrix}
\tilde A_{1}\\
\tilde A_{2}\\
\tilde F_{3}
\end{bmatrix},
\end{align}
where we have introduced the vector potentials 
\begin{align}
\label{VP}
\begin{split}
\tilde A_{k}&=\frac{1}{k_1^2+k_2^2-\omega^2\varepsilon_{\rm b}\mu_{\rm b}}\tilde
J_k^{\rm ind},\quad k=1,2,3;
\\
\tilde F_{m}&=\frac{1}{k_1^2+k_2^2-\omega^2\varepsilon_{\rm b}\mu_{\rm b}}\tilde
K_m^{\rm ind},\quad m=1,2,3.
\end{split}
\end{align}
Transforming \eqref{ks_TE_2} back to the $(\bx,\omega)$ domain and recognizing the partial
derivatives as $k_n\rightarrow i\partial_n,\,n=1,2$, we get
\begin{align}
\label{TE_sc}
\begin{bmatrix}
E_{1}^{\rm sc}\\
E_{2}^{\rm sc}\\
H_{3}^{\rm sc}
\end{bmatrix}=
\begin{bmatrix}
-\frac{\partial_1^2}{i\omega\varepsilon_{\rm b}}-\frac{\omega^2\varepsilon_{\rm
b}\mu_{\rm b}}{i\omega\varepsilon_{\rm b}} &
-\frac{\partial_1\partial_2}{i\omega\varepsilon_{\rm b}} & -\partial_2\\
-\frac{\partial_1\partial_2}{i\omega\varepsilon_{\rm b}} &
-\frac{\partial_2^2}{i\omega\varepsilon_{\rm b}}-\frac{\omega^2\varepsilon_{\rm
b}\mu_{\rm b}}{i\omega\varepsilon_{\rm b}} & \partial_1 \\
-\partial_2 & \partial_1 & i\omega\varepsilon_{\rm b}
\end{bmatrix}
\begin{bmatrix}
A_{1}\\
A_{2}\\
F_{3}
\end{bmatrix}.
\end{align}
The $(\bx,\omega)$--domain vector potentials are the spatial convolutions of the
induced currents with the scalar Green's function,
namely
\begin{align}
\label{VP_x}
\begin{split}
A_{k}(\bx,\omega)&=\int_{\bx'\in\mathbb R^2}g(\bx-\bx',\omega)J_k^{\rm
ind}(\bx,\omega),\quad k=1,2,3;
\\
F_{m}(\bx,\omega)&=\int_{\bx'\in\mathbb R^2}g(\bx-\bx',\omega)K_m^{\rm
ind}(\bx,\omega),\quad m=1,2,3.
\end{split}
\end{align}
The Green's function $g(\bx,\omega)$ is the two-dimensional inverse Fourier transform of
the expression $1/(k_1^2+k_2^2-\omega^2\varepsilon_{\rm b}\mu_{\rm b})$
appearing in \eqref{VP}. It is easy to verify that this Green's function satisfies the 
non-homogeneous Helmholtz equation with a line current (two-dimensional Dirac's delta function) 
as a source term, located at $\bx'$ position. Since the scattered fields are supposed
to satisfy the so-called radiation boundary condition (i.e. outgoing waves decaying at infinity),
out of the two possible solutions of the said Helmholtz equation one choses
\begin{align}
\label{g}
g(\bx-\bx',\omega)=\frac{i}{4}H_0^{(1)}(k_{\rm b}\vert\bx-\bx'\vert).
\end{align}
The other possible solution has the form of the Hankel function of the second kind, 
and is chosen when a different time convention is used, i.e. for the time-dependence 
of the form $e^{i\omega t}$.
\par
Substituting the induced currents from \eqref{eq:ContrastCurrents} and
expressing the scattered fields as $\bE^{\rm sc}=\bE-\bE^{\rm in}$ and $\bH^{\rm
sc}=\bH-\bH^{\rm in}$, we arrive at the following integro-differential
equations with the total fields as the fundamental unknown:
\begin{align}
\label{eq:MatrixSolutionTE}
\begin{bmatrix}
E_{1}^{\rm in}\\
E_{2}^{\rm in}\\
H_{3}^{\rm in}
\end{bmatrix}
&=
\begin{bmatrix}
E_{1}\\
E_{2}\\
H_{3}
\end{bmatrix}
-
\begin{bmatrix}
k_{\rm b}^{2}+\partial_{1}^{2} & \partial_{1}\partial_{2} & 
-i\omega\mu_{\rm b}(-\partial_{2})\\
\partial_{2}\partial_{1} & k_{\rm b}^{2}+\partial_{2}^{2} &
-i\omega\mu_{\rm b}\partial_{1}\\
-i\omega\varepsilon_{\rm b}(-\partial_{2}) & -i\omega\varepsilon_{\rm
b}\partial_{1} &
k_{\rm b}^{2}
\end{bmatrix}
\begin{bmatrix}
g*(\chi_{{\rm e}}E_{1})\\
g*(\chi_{{\rm e}}E_{2})\\
g*(\chi_{{\rm m}}H_{3})
\end{bmatrix},
\end{align}
where 
\begin{align}
\label{eq:ElectricContrast}
\chi_{\rm {\rm e}}(\bx,\omega)&=\frac{\varepsilon(\bx,\omega)}{\varepsilon_{\rm
b}}-1,
\\
\chi_{\rm {\rm m}}(\bx,\omega)&=\frac{\mu(\bx,\omega)}{\mu_{\rm b}}-1,
\end{align}
are the normalized electric and magnetic contrast functions, respectively, and
the star $(\ast)$ denotes the 2D convolution, for example
\begin{align}
\label{eq:Convolution}
g*(\chi_{{\rm e}}E_{1})= \int_{\bx'\in{\mathbb R}^{2}}\,
g(\bx-\bx',\omega)\chi_{{\rm e}}(\bx')E_{1}(\bx'){\rm d}\bx'.
\end{align}

\par
The system (\ref{eq:MatrixSolutionTE}) contains nine scalar integro-differential operators.
If the partial derivatives are carried out, then we arrive at nine ``pure'' integral operators
whose kernels have different degrees of singularity. The weakly singular kernels result in
compact operators, whereas the strongly singular kernels require some extra caution.
The presence of the second-order partial derivatives in the upper left $(2\times 2)$-corner of the 
derivative matrix in (\ref{eq:MatrixSolutionTE}) indicates that we have four scalar strongly 
singular kernels in the TE case. The corresponding operators are called singular integral operators.

To arrive at the standard form of the singular integral operator we first 
separate the domain of integration $D$ in two sub-domains as $D=[D\setminus D(\epsilon)]\cup
D(\epsilon)$ where $D(\epsilon)$ is a circular area around $\bx$ with
the radius $\epsilon$. Let $u_r$ be either of the electric field
components $E_1$ or $E_2$, then the product between the operator matrix and the
convolution vector in \eqref{eq:MatrixSolutionTE} can be written as
\begin{align}
\label{sD}
I_1=&\lim_{\epsilon\rightarrow0}\int_{\bx'\in[D\setminus
D(\epsilon)]}\partial_k\partial_r g(\bx-\bx')\chi_{\rm e}(\bx')u_r(\bx'){\rm
d}\bx'\notag\\
&+\lim_{\epsilon\rightarrow0}\partial_k\int_{\bx'\in
D(\epsilon)}\partial_rg(\bx-\bx')\chi_{\rm e}(\bx')u_r(\bx'){\rm d}\bx',\quad
k,r=1,2.
\end{align}
Here and in what follows we save some space by not showing the parametric 
dependence on $\omega$.
The first term in \eqref{sD} is recognized as the principal value, while the
second integral will be denoted by $I_2$. This second term incorporates the
Green's function of \eqref{g}. Utilizing asymptotic expansions for small
arguments for the zero order Bessel and Neumann functions \cite{abramowitz},
namely
\begin{align}
\label{J0}
J_0(z)&=1-\frac{z^2}{4}+O(z^4),\\
\label{N0}
N_0(z)&=\frac{2}{\pi}\gamma J_0(z)+\frac{2}{\pi}J_0(z)\ln\frac{z}{2}+O(z^2),
\end{align}
we replace the Hankel function in Green's function by its asymptotic
expansion based on \eqref{J0} and \eqref{N0}, arriving at
\begin{align}
\label{g_exp}
g(\bx-\bx')&=-\frac{\gamma}{2\pi}+\frac{i}{4}-\frac{1}{2\pi}\ln\frac{k_{\rm
b}\vert\bx-\bx'\vert}{2}+O(\vert\bx-\bx'\vert^2)\notag\\
&=\,g_0(\bx-\bx')+g_1(\bx-\bx'),
\end{align}
where we have defined
\begin{align}
g_0(\bx-\bx')=g(\bx-\bx')-\left[-\frac{1}{2\pi}\ln\frac{k_{\rm
b}\vert\bx-\bx'\vert}{2}\right],
\end{align}
and
\begin{align}
g_1(\bx-\bx')=-\frac{1}{2\pi}\ln\frac{k_{\rm b}\vert\bx-\bx'\vert}{2}.
\end{align}
The function $g_0(\bx-\bx')$ is not singular and its contribution (even after differentiation) 
in $I_2$ will vanish in the limit $\epsilon\rightarrow0$. The nonzero contribution to $I_2$
comes from $g_1(\bx-\bx')$. Taking the first partial derivative $\partial_r$
of $g_1(\bx-\bx')$, we get
\begin{align}
\label{g1_p}
\partial_r g_1(\bx-\bx')=\frac{1}{2\pi}\frac{-\Theta_r}{\vert\bx-\bx'\vert},
\end{align}
where $\Theta_r=(x_r-x_r')/|\bx-\bx'|$. A singularity of order one has appeared,
which is still a weak singularity (the order of the singularity is smaller than the dimension of
the manifold, which is two-dimensional in the present problem). 
Applying the second partial derivative $\partial_k$ we obtain a strong 
second-order singularity:
\begin{align}
\label{g1_pp}
\partial_k\partial_r
g_1(\bx-\bx')=\frac{1}{2\pi}\partial_k\frac{-\Theta_r}{\vert\bx-\bx'\vert}=\frac
{1}{2\pi}\frac{1}{\vert\bx-\bx'\vert^2}\left(2\Theta_k\Theta_r-\delta_{kr}
\right).
\end{align}
Thus, omitting the $g_{0}$ part (since it disappears in the limit) 
we re-write the second term of \eqref{sD} as
\begin{align}
\label{I3}
&\lim_{\epsilon\rightarrow0}\partial_k\int_{\bx'\in
D(\epsilon)}\partial_rg(\bx-\bx')\chi_{\rm e}(\bx')u_r(\bx'){\rm d}\bx'\notag\\
&=\lim_{\epsilon\rightarrow0}\int_{\bx'\in
D(\epsilon)}\frac{1}{2\pi}\partial_k\frac{-\Theta_r}{\vert\bx-\bx'\vert}\chi_{
\rm e}(\bx')u_r(\bx'){\rm d}\bx'.
\end{align}
By adding and subtracting $\chi_{\rm e}(\bx)u_r(\bx)$ inside the integral of
\eqref{I3}, we get
\begin{align}
\label{I4}
&\lim_{\epsilon\rightarrow0}\partial_k\int_{\bx'\in
D(\epsilon)}\partial_rg(\bx-\bx')\chi_{\rm e}(\bx')u_r(\bx'){\rm d}\bx'\notag\\
&=\frac{1}{2\pi}\lim_{\epsilon\rightarrow0}\int_{\bx'\in
D(\epsilon)}\partial_k\frac{-\Theta_r}{\vert\bx-\bx'\vert}\left[\chi_{\rm
e}(\bx')u_r(\bx')-\chi_{\rm e}(\bx)u_r(\bx)\right]{\rm d}\bx'\notag\\
&+\frac{1}{2\pi}\lim_{\epsilon\rightarrow0}\chi_{\rm
e}(\bx)u_r(\bx)\int_{\bx'\in
D(\epsilon)}\partial_k\frac{-\Theta_r}{\vert\bx-\bx'\vert}{\rm d}\bx'.
\end{align}
Assuming the H{\"o}lder continuity of the function $\chi_{\rm e}(\bx)u_r(\bx)$,
i.e., assuming that there exist $\alpha,C >0$, such that for all $\bx,\bx'\in{\mathbb R}^{2}$,
\begin{align}
 \label{eq:Holder}
 \vert\chi_{\rm e}(\bx)u_r(\bx)-\chi_{\rm e}(\bx')u_r(\bx')\vert \le C\vert \bx-\bx'\vert^{\alpha},
\end{align}
we effectively lower the order of singularity
in the first integral in the RHS of \eqref{I4}. Hence, in the limit $\epsilon\rightarrow 0$ this term is
zero. 
Interchanging $\partial_k$ with $-\partial_k'$ and applying the 2D
divergence theorem in the remaining term of \eqref{I4}, we obtain
\begin{align}
\label{I5}
&\lim_{\epsilon\rightarrow0}\partial_k\int_{\bx'\in
D(\epsilon)}\partial_rg(\bx-\bx')\chi_{\rm e}(\bx')u_r(\bx'){\rm d}\bx'\notag\\
&=\frac{1}{2\pi}\chi_{\rm
e}(\bx)u_r(\bx)\lim_{\epsilon\rightarrow0}\oint_{\bx'\in \partial
D(\epsilon)}-\nu_k\frac{-\Theta_r}{\vert\bx-\bx'\vert}{\rm d}\bx',
\end{align}
with $\nu_k=-\Theta_k$ being the projection of the normal unit vector of the
polar coordinate system on the $x_k$ axis. Finally, we arrive at the result
\begin{align}
\label{I6}
&\lim_{\epsilon\rightarrow0}\partial_k\int_{\bx'\in
D(\epsilon)}\partial_rg(\bx-\bx')\chi_{\rm e}(\bx')u_r(\bx'){\rm d}\bx'\notag\\
&=\frac{1}{2\pi}\chi_{\rm
e}(\bx)u_r(\bx)\lim_{\epsilon\rightarrow0}\oint_{\bx'\in \partial
D(\epsilon)}-\frac{\Theta_k\Theta_r}{\vert\bx-\bx'\vert}{\rm d}\bx'\notag\\
&=-\frac{1}{2\pi}\chi_{\rm
e}(\bx)u_r(\bx)\lim_{\epsilon\rightarrow0}\oint_{\bx'\in \partial
D(1)}\Theta_k\Theta_r{\rm d}\bx'\notag\\
&=-\frac{1}{2}\chi_{\rm e}(\bx)u_r(\bx)\delta_{rk},
\end{align}
where we notice that the last contour integral is over the unit circle. Hence,
\eqref{sD} can now be written as
\begin{align}
I_1=&\,p.\,v.\int_{\bx'\in D}\partial_k\partial_rg(\bx-\bx')\chi_{\rm
e}(\bx')u_r(\bx'){\rm d}\bx'\notag\\
&-\frac{1}{2}\chi_{\rm e}(\bx)u_r(\bx)\delta_{rk},\quad k,r=1,2.
\end{align}
Now, we can represent \eqref{eq:MatrixSolutionTE} in the standard form
\begin{align}
\label{eq:SingularTE}
\begin{split}
\begin{bmatrix}
E_{1}^{\rm in}\\
E_{2}^{\rm in}\\
H_{3}^{\rm in}
\end{bmatrix}
=
\begin{bmatrix}
S & 0 & 0 \\
0 & S & 0 \\
0 & 0 & I
\end{bmatrix}
\begin{bmatrix}
E_{1}\\
E_{2}\\
H_{3}
\end{bmatrix}
&
+
p.\,v.
\begin{bmatrix}
G_{11} & G_{12} & 0 \\
G_{21} & G_{22} & 0 \\
0 & 0 & 0
\end{bmatrix}*
\begin{bmatrix}
X_{{\rm e}}E_{1}\\
X_{{\rm e}}E_{2}\\
X_{{\rm m}}H_{3}
\end{bmatrix}
\\
&+
\begin{bmatrix}
K_{11} & K_{12} & K_{13} \\
K_{21} & K_{22} & K_{23} \\
K_{31} & K_{32} & K_{33}
\end{bmatrix}*
\begin{bmatrix}
X_{{\rm e}}E_{1}\\
X_{{\rm e}}E_{2}\\
X_{{\rm m}}H_{3}
\end{bmatrix},
\end{split}
\end{align}
with $S$ denoting the operator of (pointwise) multiplication with the function $s(\bx)=1+1/2\chi_{\rm e}(\bx)$, 
and $I$ -- the identity operator. 
The kernels in the principal-value operator in \eqref{eq:SingularTE} are easily recognized 
using \eqref{g1_pp} and
\eqref{eq:MatrixSolutionTE}, and are given by
\begin{align}
G_{nm}(\bx)=-\frac{1}{2\pi\vert\bx\vert^2}\left[2\Theta_n\Theta_m-\delta_{nm}
\right],\quad n,m=1,2.
\end{align}
The kernels in the compact operator in \eqref{eq:SingularTE} are easily obtained
using \eqref{pd_1}, \eqref{pd_2} and \eqref{eq:MatrixSolutionTE}, and are given
by
\begin{align}
\label{eq:KompactKernel1122}
\begin{split}
K_{nm}=&\left[\frac{1}{2\pi\vert\bx\vert^2}-\frac{ik_{\rm
b}}{4\vert\bx\vert}H_{1}^{(1)}(k_{\rm b}\vert\bx\vert)
\right]
\left[2\Theta_{n}\Theta_{m}-\delta_{nm}\right]
\\
&+\frac{ik_{\rm b}^2}{4}H_{0}^{(1)}(k_{\rm b}\vert\bx\vert)
\left[\Theta_{n}\Theta_{m}-\delta_{nm}\right], \;\;\;n,m=1,2;
\end{split}
\\
\label{eq:KompactKernel13}
\begin{split}
K_{13}=-\frac{\omega\mu_{\rm b}k_{\rm b}\Theta_{2}}{4}H_{1}^{(1)}(k_{\rm
b}\vert\bx\vert),
\end{split}
\\
\label{eq:KompactKernel31}
\begin{split}
K_{31}=-\frac{\omega\varepsilon_{\rm b}k_{\rm b}\Theta_{2}}{4}H_{1}^{(1)}(k_{\rm
b}\vert\bx\vert),
\end{split}
\\
\label{eq:KompactKernel23}
\begin{split}
K_{23}=\frac{\omega\mu_{\rm b}k_{\rm b}\Theta_{1}}{4}H_{1}^{(1)}(k_{\rm
b}\vert\bx\vert),
\end{split}
\\
\label{eq:KompactKernel32}
\begin{split}
K_{32}=\frac{\omega\varepsilon_{\rm b}k_{\rm b}\Theta_{1}}{4}H_{1}^{(1)}(k_{\rm
b}\vert\bx\vert),
\end{split}
\\
\label{eq:KompactKernel33}
\begin{split}
K_{33}=-\frac{ik_{\rm b}^{2}}{4}H_{0}^{(1)}(k_{\rm b}\vert\bx\vert).
\end{split}
\end{align}

\section{\label{sec:symbol}Derivation of the Symbol}
Using operator notation, \eqref{eq:SingularTE} can be written as
\begin{align}
\label{ie_mn}
(A\bE)(\bx)=\left(\mathbb E+\frac{1}{2}\mathbb M\right)\bE(\bx)+\mathbb A^{\rm
s}\mathbb M\bE(\bx')+\mathbb K\bE(\bx'),
\end{align}
where $\mathbb A^{\rm s}$ is a (matrix) singular integral operator, $\mathbb K$ is a
(matrix) compact integral operator, $\mathbb E$ is a (matrix) identity operator, and $\mathbb
M(\bx)=\chi_{\rm e}(\bx){\mathbb E}_{2}$ is a (matrix) multiplication operator, where the 
lower-right element of ${\mathbb E}_{2}$ is set to zero.
The singular term of \eqref{ie_mn} is of the form
\begin{align}
(A^{\rm s}\mathbb M\bE)(\bx)=\int_{\bx'\in
D}\frac{{\mathbb F}(\mathbf\Theta)}{\vert\bx-\bx'\vert^2}\chi_{\rm e}(\bx')\bE(x'){\rm
d}\bx'.
\end{align}
According to \eqref{g1_pp}, the characteristic (matrix) function ${\mathbb F}(\mathbf\Theta)$ 
has the form
\begin{align}
\label{char}
{\mathbb F}(\mathbf\Theta)=-\frac{1}{2\pi}[2\mathbb
Q(\bx-\bx')-{\mathbb I}_{2}],
\end{align}
where ${\mathbb I}_{2}$ is the $(3\times 3)$ identity matrix with the lower-right element 
set to zero, and the tensor $\mathbb Q$ is given by
\begin{align}
\label{Q}
\mathbb Q(\bx-\bx')=
\begin{bmatrix}
\frac{(x_1-x_1')^2}{\vert\bx-\bx'\vert^2} &
\frac{(x_1-x_1')(x_2-x_2')}{\vert\bx-\bx'\vert^2} & 0\\
\frac{(x_1-x_1')(x_2-x_2')}{\vert\bx-\bx'\vert^2} &
\frac{(x_2-x_2')^2}{\vert\bx-\bx'\vert^2} & 0\\
0 & 0 & 0
\end{bmatrix}.
\end{align}
For $\chi_{\rm e}(\bx)$ H{\"o}lder-continuous on ${\mathbb R}^{2}$, 
the symbol of the compound singular integral operator of \eqref{ie_mn} 
can be computed as (see \cite{Mikhlin}):
\begin{align}
\label{symbol_b}
{\rm Smb}(\mathbb A)=\mathbb I+\frac{1}{2}\chi_{\rm
e}(\bx){\mathbb I}_{2}+{\rm Smb}(\mathbb A^{\rm s})\chi_{\rm e}(\bx),
\end{align}
where ${\mathbb I}$ is the ordinary $(3\times 3)$ identity matrix.
Since the characteristics ${\mathbb F}$ depends only 
on $\bx-\bx'$, the symbol of the singular
integral operator $\mathbb A^{\rm s}$ is the Fourier transform of its kernel
${\mathbb F}(\mathbf\Theta)/|\bx-\bx'|^2$ with respect to the 
variable $\by=\bx-\bx'$, i.e. it is a $\bk$-domain matrix-valued function 
$\tilde{\mathbb F}^{\rm s}(\bk)$. Computing this Fourier transform is a 
daunting task, and we shall use a shortcut proposed in \cite{Mikhlin}.

Let $A^{\rm s}$ be a single component of our matrix-valued operator 
${\mathbb A}^{\rm s}$ and let $\mathbf\Phi^{\rm s}(\bk)$ denote its symbol,
which is one of the components of the matrix-valued symbol function $\tilde{\mathbb F}^{\rm s}(\bk)$
we are trying to compute.
The symbol ${\rm Smb}(A^{\rm s})=\mathbf\Phi^{\rm s}(\bk)$ of 
a scalar simple singular integral operator can be expanded in a series of 2D spherical functions 
of order $p$, that is, in Fourier series of sines and cosines \cite{Mikhlin}, namely
\begin{align}
\label{symbol_exp}
\mathbf\Phi^{\rm
s}(\tilde{\bTheta})=\sum_{p=0}^\infty\left[\gamma_{2,p}a_p^{(1)}Y_{p,2}^{(1)}
(\tilde{\bTheta})+\gamma_{2,p}a_p^{(2)}Y_{p,2}^{(2)}(\tilde{\bTheta})\right],
\end{align}
where $\tilde{\bTheta}=\bk/|\bk|$ is the unit vector in the $\bk$-domain,
$Y_{p,2}^{(1)}=\sin(p\tilde\phi)$, and $Y_{p,2}^{(2)}=\cos(p\tilde\phi)$ is the
basis of the expansion, $a_p^{(1)}$ and $a_p^{(2)}$ are the expansion
coefficients, and $\tilde\phi$ is the directional angle of the unit vector
$\tilde{\bTheta}$. In \eqref{symbol_exp}, $\gamma_{2,p}=\pi
i^p\Gamma(p/2)/\Gamma((2+p)/2)$ \cite{Mikhlin}. Since any component $f(\bTheta)$ of 
the characteristic matrix-valued function ${\mathbb F}(\bTheta)$ given by \eqref{char} 
depend only on $\bTheta$, we can expand each of them in a Fourier series as well
\begin{align}
\label{char_exp}
f(\mathbf\Theta)=\sum_{p=1}^\infty\left[a_p^{(1)}\sin(p\phi)+a_p^{(2)}
\cos(p\phi)\right],
\end{align}
where $\phi$ is the directional angle of the unit vector $\bTheta=(\bx-\bx')/|\bx-\bx'|$. 
It was shown in \cite{Mikhlin} that the
expansion coefficients in \eqref{char_exp} are the same with those of
\eqref{symbol_exp}. Using \eqref{char} and \eqref{Q}, we see that the components
of the characteristic ${\mathbb F}$ are
\begin{align}
\label{f_m}
{\mathbb F}(\mathbf\Theta)=-\frac{1}{2\pi}[2\mathbb Q(\bx-\bx')-{\mathbb
I}_{2}]=-\frac{1}{2\pi}
\begin{bmatrix}
2\cos^2\phi-1 & 2\cos\phi\sin\phi & 0\\
2\cos\phi\sin\phi & 2\sin^2\phi-1 & 0\\
0 & 0 & 0
\end{bmatrix}.
\end{align}
To calculate the expansion coefficients, we must equate each element
$[{\mathbb F}(\mathbf\Theta)]_{kr},\,k,r=1,2$, from \eqref{f_m} with the series from
\eqref{char_exp}. For example, for the element $[{\mathbb F}(\mathbf\Theta)]_{11}$ we have
\begin{align}
-\frac{1}{2\pi}[2\cos^2\phi-1]=-\frac{1}{2\pi}\cos(2\phi)=\sum_{p=1}^\infty\left
[a_p^{(1)}\sin(p\phi)+a_p^{(2)}\cos(p\phi)\right].
\end{align}
From this expansion it is obvious that the only nonzero expansion
coefficient is $a_2^{(2)}=-1/(2\pi)$, while the rest are all zero (i.e.
$a_p^{(1)}=0$ $\forall p$ and $a_p^{(2)}=0$ $\forall p\neq2$). The same
procedure is followed for the rest of the components in \eqref{f_m}. Then, we substitute
the known expansion coefficients in \eqref{symbol_exp} to get the elements
$[\tilde{\mathbb F}^{\rm s}(\tilde{\bTheta})]_{kr},\,k,r=1,2$. Continuing our example, the element
$[\tilde{\mathbb F}^{\rm s}(\tilde{\bTheta})]_{11}$ is obtained through the substitutions
$\gamma_{2,2}=\pi i^2\Gamma(1)/\Gamma(2)=-\pi$,
$Y_{2,2}^{(2)}=\cos(2\tilde\phi)$ and $a_2^{(2)}=-1/(2\pi)$, and therefore
$[\tilde{\mathbb F}^{\rm s}(\tilde{\bTheta})]_{11}=1/2\cos(2\tilde\phi)=\cos^2(\tilde\phi)-1/2$.
Following the same procedure for all components, we finally get
\begin{align}
\label{symbol_s}
{\rm Smb}(\mathbb A^{\rm s})=\tilde{\mathbb F}^{\rm s}(\tilde{\bTheta})=
\begin{bmatrix}
\cos^2\tilde\phi-1/2 & \sin\tilde\phi\cos\tilde\phi & 0 \\
\sin\tilde\phi\cos\tilde\phi & \sin^2\tilde\phi-1/2 & 0 \\
0 & 0 & 0
\end{bmatrix}
={\mathbb Q}(\bk)-\frac{1}{2}{\mathbb I}_{2}.
\end{align}
Substituting the result of \eqref{symbol_s} back in \eqref{symbol_b}, we get the
symbol of the complete operator as the following $(3\times 3)$ matrix-valued function:
\begin{align}
\label{symbol}
{\rm Smb}(\mathbb A)(\bx,\bk)=\mathbb I+\chi_{\rm e}(\bx)\mathbb Q(\bk).
\end{align}

\section{Derivation of the algebraic system matrix}
Equation \eqref{eq:SingularTE} defines an algebraic system $Au=b$ for the
numerical evaluation of the fields. To derive the matrix elements, its more
convenient to use the equivalent integral form of \eqref{eq:SingularTE}. 

The application of the matrix operator of \eqref{eq:MatrixSolutionTE} --the
matrix that contains the partial derivatives-- on Green's functions, give us the
Green's tensor. In order to get the equivalent integral form of
\eqref{eq:SingularTE}, we need to split the Green's tensor in two parts. The
first part, denoted by $\mathbb A(\bx-\bx')$, corresponds to the second order or
mixed derivatives only, i.e.
\begin{align}
\mathbb A(\bx-\bx')=\left(k_{\rm b}^2
\begin{bmatrix}
1 & 0 & 0\\
0 & 1 & 0\\
0 & 0 & 1
\end{bmatrix}
+
\begin{bmatrix}
\partial_1^2 & \partial_1\partial_2 & 0\\
\partial_2\partial_1 & \partial_2^2 & 0\\
0 & 0 & 0
\end{bmatrix}
\right)g(\bx-\bx').
\end{align}
The second part, denoted by $\mathbb B(\bx-\bx')$, corresponds to the first
order derivatives only, i.e.
\begin{align}
\mathbb B(\bx-\bx')=
\begin{bmatrix}
0 & 0 & \partial_2\\
0 & 0 & -\partial_1\\
-\partial_2 & \partial_1 & 0
\end{bmatrix}
g(\bx-\bx').
\end{align}
To obtain  the explicit relations for $\mathbb A(\bx-\bx')$ and $\mathbb
B(\bx-\bx')$, we first calculate the first order derivative which is
\begin{align}
\label{pd_1}
\partial_rg(\bx-\bx')&=-\frac{i}{4}k_{\rm
b}\frac{x_r-x_r'}{\vert\bx-\bx'\vert}H_1^{(1)}(k_{\rm
b}\vert\bx-\bx'\vert),\quad r=1,2.
\end{align}
Then, the mixed derivative is
\begin{align}
\label{pd_2}
\partial_k\partial_rg(\bx-\bx')=&\frac{i}{4}\bigg[k_{\rm
b}\frac{1}{\vert\bx-\bx'\vert}\left(2\frac{x_r-x_r'}{\vert\bx-\bx'\vert}\frac{
x_k-x_k'}{\vert\bx-\bx'\vert}-\delta_{kr}\right)H_1^{(1)}(k_{\rm
b}\vert\bx-\bx'\vert)\notag\\
&-k_{\rm
b}^2\frac{x_r-x_r'}{\vert\bx-\bx'\vert}\frac{x_k-x_k'}{\vert\bx-\bx'\vert}H_0^{
(1)}(k_{\rm b}\vert\bx-\bx'\vert)\bigg],\quad k,r=1,2;
\end{align}
where $\delta_{kr}$ is the Kronecker's delta. So, the $\mathbb A$ tensor is
given by
\begin{align}
\label{G1}
\mathbb A(\bx-\bx')=&\frac{i}{4}k_{\rm
b}\frac{1}{\vert\bx-\bx'\vert}H_1^{(1)}(k_{\rm b}\vert\bx-\bx'\vert)[2\mathbb
Q(\bx-\bx')-\mathbb I_2]\notag\\
&-\frac{i}{4}k_{\rm b}^2H_0^{(1)}(k_{\rm b}\vert\bx-\bx'\vert)[\mathbb
Q(\bx-\bx')-\mathbb I_2].
\end{align}
The tensor $\mathbb Q$ was introduced in \eqref{Q} while $\mathbb I_2$
is, as explained in \sec{sec:symbol}, the $(3\times 3)$ identity matrix with the lower-right element 
set to zero. The $\mathbb B$ tensor is obtained with the use of
\eqref{pd_1}
\begin{align}
\label{G2}
\mathbb B(\bx-\bx')=-\frac{i}{4}k_{\rm b}H_1^{(1)}(k_{\rm
b}\vert\bx-\bx'\vert)\mathbf\Theta(\bx-\bx')\times,
\end{align}
where we have now introduced the $\mathbf\Theta\times$ tensor given by
\begin{align}
\label{Theta}
\mathbf\Theta(\bx-\bx')\times=
\begin{bmatrix}
0 & 0 & \frac{x_2-x_2'}{\vert\bx-\bx'\vert}\\
0 & 0 & -\frac{x_1-x_1'}{\vert\bx-\bx'\vert}\\
-\frac{x_2-x_2'}{\vert\bx-\bx'\vert} & \frac{x_1-x_1'}{\vert\bx-\bx'\vert} & 0
\end{bmatrix}.
\end{align}
\par
After this splitting, we can easily express \eqref{eq:SingularTE} in the
following integral form
\begin{align}
\label{eq:SingularTE_i_E}
\bE^{\rm in}(\bx)=&\left[1+\frac{1}{2}\chi_{\rm e}(\bx)\right]\bE(\bx)\notag\\
&-\,p.v.\int_{\bx'\in D}\mathbb A(\bx-\bx')\chi_{\rm e}(\bx')\bE(\bx'){\rm
d}\bx'\notag\\
&-i\omega\mu_{\rm b}\int_{\bx'\in D}\mathbb B(\bx-\bx')\chi_{\rm
m}(\bx')\bH(\bx'){\rm d}\bx',\\
\label{eq:SingularTE_i_H}
\bH^{\rm in}(\bx)=&\bH(\bx)-\int_{\bx'\in D}\mathbb A(\bx-\bx')\chi_{\rm
m}(\bx')\bH(\bx'){\rm d}\bx'\notag\\
&+i\omega\varepsilon_{\rm b}\int_{\bx'\in D}\mathbb B(\bx-\bx')\chi_{\rm
e}(\bx')\bE(\bx'){\rm d}\bx',
\end{align}
where the vectors $\bE(\bx)=[E_1(\bx),E_2(\bx),0]^T$ and $\bH(\bx)=[0,0,H_3(\bx)]^T$. Same applies for $\bE^{\rm in}(\bx)$ and $\bH^{\rm in}(\bx)$.
\par
Equation \eqref{eq:SingularTE_i_E} can be rewritten as
\begin{align}
\label{TE_E}
\bE^{\rm in}(\bx)=&\left[1+\frac{1}{2}\chi_{\rm e}(\bx)\right]\bE(\bx)\notag\\
&-p.v.\int_{\bx'\in D_{\rm s}}\mathbb A(\bx-\bx')\chi_{\rm e}(\bx')\bE(\bx'){\rm
d}\bx'\notag\\
&-\int_{\bx'\in D\setminus D_{\rm s}}\mathbb A(\bx-\bx')\chi_{\rm
e}(\bx')\bE(\bx'){\rm d}\bx'\notag\\
&-i\omega\mu_{\rm b}\int_{\bx'\in D_{\rm s}}\mathbb B(\bx-\bx')\chi_{\rm
m}(\bx')\bH(\bx'){\rm d}\bx'\notag\\
&-i\omega\mu_{\rm b}\int_{\bx'\in D\setminus D_{\rm s}}\mathbb
B(\bx-\bx')\chi_{\rm m}(\bx')\bH(\bx'){\rm d}\bx',
\end{align}
where we have separated the domain $D$ into a domain $D\setminus D_{\rm s}$
``free'' of singularity, and a domain $D_{\rm s}$ which encloses the position
vector $\bx$ and hence the singularity. In this way we have three integrals in
the usual sense and one in the sense of principal value. Applying a simple
collocation technique with the mid--point rule, the usual sense integrals over
the domain $D\setminus D_{\rm s}$ would be given by
\begin{align}
\label{R1}
&-\int_{\bx'\in D\setminus D_{\rm s}}\mathbb A(\bx-\bx')\chi_{\rm
e}(\bx')\bE(\bx'){\rm d}\bx'\notag\\
&\approx-\mathop{\sum_{m=1}^{N-1}}_{m\neq n}\mathbb A(\bx_n-\bx_m)\chi_{\rm
e}(\bx_m)\bE(\bx_m)S_m,\quad n=1,2,\ldots,N;\\
\label{R2}
&-i\omega\mu_{\rm b}\int_{\bx'\in D\setminus D_{\rm s}}\mathbb
B(\bx-\bx')\chi_{\rm m}(\bx')\bH(\bx'){\rm d}\bx'\notag\\
&\approx-i\omega\mu_{\rm b}\mathop{\sum_{m=1}^{N-1}}_{m\neq n}\mathbb
B(\bx_n-\bx_m)\chi_{\rm m}(\bx_m)\bH(\bx_m)S_m,\quad n=1,2,\ldots,N;
\end{align}
where $S_m=h^2$ is the surface of each elementary cell $D_m$ of side $h$ in the
computational domain.
\par
We now proceed to calculate the principal value integral in \eqref{TE_E}. We
have
\begin{align}
\label{pv_c}
&p.v.\int_{\bx'\in D_{\rm s}}\mathbb A(\bx-\bx')\chi_{\rm e}(\bx')\bE(\bx'){\rm
d}\bx'\notag\\
&=\lim_{\epsilon\rightarrow0}\int_{\bx'\in D_{\rm s}\setminus
D(\epsilon)}\mathbb A(\bx-\bx')\chi_{\rm e}(\bx')\bE(\bx'){\rm d}\bx'\notag\\
&\approx\chi_{\rm e}(\bx_n)\bE(\bx_n)\lim_{\epsilon\rightarrow0}\int_{\bx'\in
D_{\rm s}\setminus D(\epsilon)}\mathbb A(\bx-\bx'){\rm d}\bx',
\end{align}
where the domain of integration in the last integral is over the singular cell
but with an exception of a small circular neighborhood having radius
$\epsilon$. We further define the tensor $\mathbb L$ as
\begin{align}
\label{L}
\mathbb L=\lim_{\epsilon\rightarrow0}\int_{\bx'\in D_{\rm s}\setminus
D(\epsilon)}\mathbb A(\bx-\bx'){\rm d}\bx'.
\end{align}
Now we recall the form that $\mathbb A(\bx-\bx')$ has from \eqref{G1}. It is
convenient to rewrite \eqref{G1} in terms of $\mathbb Q(\bx-\bx')$ and
$\mathbb I_2$, namely
\begin{align}
\label{G1_QI}
\mathbb A(\bx-\bx')=&\frac{i}{4}\left[2k_{\rm
b}\frac{1}{\vert\bx-\bx'\vert}H_1^{(1)}(k_{\rm b}\vert\bx-\bx'\vert)-k_{\rm
b}^2H_0^{(1)}(k_{\rm b}\vert\bx-\bx'\vert)\right]\mathbb Q(\bx-\bx')\notag\\
&+\frac{i}{4}\left[k_{\rm b}^2H_0^{(1)}(k_{\rm b}\vert\bx-\bx'\vert)-k_{\rm
b}\frac{1}{\vert\bx-\bx'\vert}H_1^{(1)}(k_{\rm
b}\vert\bx-\bx'\vert)\right]\mathbb I_2.
\end{align}
To compute \eqref{L}, we transform the singular cell $D_{\rm s}$ from
rectangular shape to a circular disk having center at $\bx$ and radius $a_n$.
The transformed circular cell with the original one have the same surface $S=\pi
a_n^2=h^2$. Then, \eqref{G1_QI} is transformed into polar coordinates as
\begin{align}
\label{G1_QI_p}
\mathbb A(\rho,\phi)=&\frac{i}{4}\left[2k_{\rm b}\frac{1}{\rho}H_1^{(1)}(k_{\rm
b}\rho)-k_{\rm b}^2H_0^{(1)}(k_{\rm b}\rho)\right]\mathbb Q(\phi)\notag\\
&+\frac{i}{4}\left[k_{\rm b}^2H_0^{(1)}(k_{\rm b}\rho)-k_{\rm
b}\frac{1}{\rho}H_1^{(1)}(k_{\rm b}\rho)\right]\mathbb I_2.
\end{align}
Now, we calculate each element of the tensor $\mathbb L$ from
\begin{align}
\label{Lpq}
L_{pq}=&\frac{i}{4}\bigg\{\lim_{\epsilon\rightarrow0}\int_{\rho=\epsilon}^{a_m}
\left[2k_{\rm b}\frac{1}{\rho}H_1^{(1)}(k_{\rm b}\rho)-k_{\rm
b}^2H_0^{(1)}(k_{\rm b}\rho)\right]\rho{\rm
d}\rho\int_{\phi=0}^{2\pi}Q_{pq}(\phi){\rm d}\phi\notag\\
&+\int_{\rho=\epsilon}^{a_m}\left[k_{\rm b}^2H_0^{(1)}(k_{\rm b}\rho)-k_{\rm
b}\frac{1}{\rho}H_1^{(1)}(k_{\rm b}\rho)\right]\rho{\rm
d}\rho\int_{\phi=0}^{2\pi}\delta_{pq}{\rm d}\phi\bigg\},\quad p,q=1,2.
\end{align}
The integration in \eqref{Lpq} is carried out only three times since the tensor
$\mathbb Q(\bx-\bx')$ is symmetric. The angular integral with integrand function
$Q_{pq}(\phi)$ gives $\pi\delta_{pq}$, while the other with the integrand
function $\delta_{pq}$ gives $2\pi\delta_{pq}$. The radial integrals are
evaluated by changing variables $z=k_{\rm b}\rho$ and using the well known
integral formulas for Bessel functions \cite{watson}
\begin{align}
\int^z\zeta^{n+1}Z_n(\zeta){\rm d}\zeta=z^{n+1}Z_{n+1}(z),\\
\int^z\zeta^{-n+1}Z_n(\zeta){\rm d}\zeta=-z^{-n+1}Z_{n-1}(z).
\end{align}
The result is
\begin{align}
L_{pq}=\delta_{pq}\left\{\frac{i\pi a_n}{4}k_{\rm b}H_1^{(1)}(k_{\rm
b}a_n)-\frac{i\pi}{4}k_{\rm b}\lim_{\epsilon\rightarrow0}\left[\epsilon
H_1^{(1)}(k_{\rm b}\epsilon)\right]\right\},\quad p,q=1,2.
\end{align}
The involved limit is easily calculated giving $-i2/(\pi k_{\rm b})$. Therefore,
substituting the results for the tensor $\mathbb L$ back to \eqref{pv_c}
\begin{align}
\label{R3}
&p.v.\int_{\bx'\in D_{\rm s}}\mathbb A(\bx-\bx')\chi_{\rm e}(\bx')\bE(\bx'){\rm
d}\bx'\notag\\
&\approx\chi_{\rm e}(\bx_n)\left[-\frac{1}{2}+\frac{i\pi a_n}{4}k_{\rm
b}H_1^{(1)}(k_{\rm b}a_n)\right]\bE(\bx_n),\quad\bx_n\in D_{\rm s}.
\end{align}
\par
What remains is the calculation of the integral involving $\mathbb B(\bx-\bx')$
over the domain $D_{\rm s}$ in \eqref{TE_E}. We have
\begin{align}
&-i\omega\mu_{\rm b}\int_{\bx'\in D_{\rm s}}\mathbb B(\bx-\bx')\chi_{\rm
m}(\bx')\bH(\bx'){\rm d}\bx'\notag\\
&\approx-i\omega\mu_{\rm b}\chi_{\rm m}(\bx_n)\bH(\bx_n)\int_{\bx'\in D_{\rm
s}}\mathbb B(\bx-\bx'){\rm d}\bx'.
\end{align}
Utilizing \eqref{G2} and \eqref{Theta}, and expressing $\mathbb B(\bx-\bx')$ in
polar coordinates, in the same manner as we did for the $\mathbb L$ tensor
before, we arrive at an angular integration from $0$ to $2\pi$ for each
component $\Theta_{pq}$ of $\mathbf\Theta(\bx-\bx')\times$ tensor. This
integration is zero for each component $\Theta_{pq}$, and hence this term does
not contribute.
\par
Putting together the results of \eqref{R1}, \eqref{R2} and \eqref{R3}, we arrive
at the linear equations
\begin{align}
\label{S1}
&\left\{1-\chi_{\rm e}(\bx_n)\left[\frac{i\pi k_{\rm
b}h}{4\sqrt{\pi}}H_1^{(1)}(k_{\rm
b}h/\sqrt{\pi})-1\right]\right\}\bE(\bx_n)\notag\\
&-h^2\mathop{\sum_{m=1}^{N-1}}_{m\neq n}\mathbb A(\bx_n-\bx_m)\chi_{\rm
e}(\bx_m)\bE(\bx_m)\notag\\
&-i\omega\mu_{\rm b}h^2\mathop{\sum_{m=1}^{N-1}}_{m\neq n}\mathbb
B(\bx_n-\bx_m)\chi_{\rm m}(\bx_m)\bH(\bx_m)=\bE^{\rm in}(\bx_n),\quad
n=1,2,\ldots,N.
\end{align}
Following the same analysis for \eqref{eq:SingularTE_i_H}, we arrive at the
linear equations
\begin{align}
\label{S2}
&\left\{1-\chi_{\rm m}(\bx_n)\left[\frac{i\pi k_{\rm
b}h}{2\sqrt{\pi}}H_1^{(1)}(k_{\rm
b}h/\sqrt{\pi})-1\right]\right\}\bH(\bx_n)\notag\\
&-h^2\mathop{\sum_{m=1}^{N-1}}_{m\neq n}\mathbb A(\bx_n-\bx_m)\chi_{\rm
m}(\bx_m)\bH(\bx_m)\notag\\
&+i\omega\varepsilon_{\rm b}h^2\mathop{\sum_{m=1}^{N-1}}_{m\neq n}\mathbb
B(\bx_n-\bx_m)\chi_{\rm e}(\bx_m)\bE(\bx_m)=\bH^{\rm in}(\bx_n),\quad
n=1,2,\ldots,N.
\end{align}
Equations \eqref{S1} and \eqref{S2} compose the linear system for the evaluation
of the unknown total fields, with the following structure
\begin{align}
\label{sm}
\begin{bmatrix}
A_{11} & A_{12} & A_{13}\\
A_{21} & A_{22} & A_{23}\\
A_{31} & A_{32} & A_{33}
\end{bmatrix}
\begin{bmatrix}
u_1\\u_2\\u_3
\end{bmatrix}
=
\begin{bmatrix}
b_1\\b_2\\b_3
\end{bmatrix}.
\end{align}
In \eqref{sm}, $[u_1,u_2,u_3]^T=[E_1(\bx_n),E_2(\bx_n),H_3(\bx_n)]^T$,
$n=1,2,\ldots,N$, are the unknown total field components on the grid, while
$[b_1,b_2,b_3]^T=[E_1^{\rm in}(\bx_n),E_2^{\rm in}(\bx_n),H_3^{\rm
in}(\bx_n)]^T$ contains the grid values of the incident field. The elements of
the system matrix in \eqref{sm} are now easily recognized with the use of
\eqref{S1}, \eqref{S2}, \eqref{G1}, \eqref{G2} and \eqref{g}. The results are
\begin{align}
\label{A11}
[A_{\ell\ell}]_{nm}&=-k_{\rm b}^2h^2\chi_{\rm
e}(\bx_m)\bigg\{\left[\frac{i}{2k_{\rm b}r_{nm}}H_1^{(1)}(k_{\rm
b}r_{nm})-\frac{i}{4}H_0^{(1)}(k_{\rm
b}r_{nm})\right]\theta_{\ell,nm}\theta_{\ell,nm}\notag\\
&+\left[\frac{i}{4}H_0^{(1)}(k_{\rm b}r_{nm})-\frac{i}{4k_{\rm
b}r_{nm}}H_1^{(1)}(k_{\rm b}r_{nm})\right]\delta_{\ell\ell}\bigg\},\quad
\ell=1,2,\quad m\neq n;
\end{align}
\begin{align}
[A_{\ell\ell}]_{nn}=&1+\left[1-\frac{i\pi k_{\rm
b}h}{4\sqrt{\pi}}H_1^{(1)}(k_{\rm b}h/\sqrt{\pi})\right]\chi_{\rm
e}(\bx_n),\quad \ell=1,2;
\end{align}
\begin{align}
[A_{\ell q}]_{nm}=&-k_{\rm b}^2h^2\chi_{\rm e}(\bx_m)\left[\frac{i}{2k_{\rm
b}r_{nm}}H_1^{(1)}(k_{\rm b}r_{nm})-\frac{i}{4}H_0^{(1)}(k_{\rm
b}r_{nm})\right]\notag\\
&\times\theta_{\ell,nm}\theta_{q,nm}(1-\delta_{nm}),\quad\ell,q=1,2,
\quad\ell\neq q;
\end{align}
\begin{align}
[A_{13}]_{nm}=&i\omega\mu_{\rm b}h^2\chi_{\rm m}(\bx_m)\frac{ik_{\rm
b}}{4}H_1^{(1)}(k_{\rm b}r_{nm})\theta_{2,nm}(1-\delta_{nm});
\end{align}
\begin{align}
[A_{32}]_{nm}=&-i\omega\varepsilon_{\rm b}h^2\chi_{\rm e}(\bx_m)\frac{ik_{\rm
b}}{4}H_1^{(1)}(k_{\rm b}r_{nm})\theta_{1,nm}(1-\delta_{nm});
\end{align}
\begin{align}
[A_{31}]_{nm}=&i\omega\varepsilon_{\rm b}h^2\chi_{\rm e}(\bx_m)\frac{ik_{\rm
b}}{4}H_1^{(1)}(k_{\rm b}r_{nm})\theta_{2,nm}(1-\delta_{nm});
\end{align}
\begin{align}
[A_{23}]_{nm}=&-i\omega\mu_{\rm b}h^2\chi_{\rm m}(\bx_m)\frac{ik_{\rm
b}}{4}H_1^{(1)}(k_{\rm b}r_{nm})\theta_{1,nm}(1-\delta_{nm});
\end{align}
\begin{align}
[A_{33}]_{nm}=-k_{\rm b}^2h^2\chi_{\rm m}(\bx_m)\frac{i}{4}H_0^{(1)}(k_{\rm
b}r_{nm}),\quad m\neq n;
\end{align}
\begin{align}
\label{A33}
[A_{33}]_{nn}=1+\left[1-\frac{i\pi k_{\rm b}h}{2\sqrt{\pi}}H_1^{(1)}(k_{\rm
b}h/\sqrt{\pi})\right]\chi_{\rm m}(\bx_n).
\end{align}
In \eqref{A11}--\eqref{A33} we have defined $r_{nm}=\vert\bx_n-\bx_m\vert$,
$\theta_{\ell,nm}=(x_{\ell,n}-x_{\ell,m})/r_{nm}$, $\ell=1,2$, $n,m=1,\ldots,N$,
and $x_{\ell,n}$ denotes the Cartesian component of the 2D position vector
$\bx_n$, pointing at the $n$th node of the grid.

\bibliographystyle{siam}
\bibliography{cleanbib}

\end{document}